\documentclass[reqno,makeidx,11pt]{amsart}
\usepackage{amssymb,amsfonts,amsmath,graphicx,
}

\def\C{{\mathbb{C}}}

\def\F{{\mathbb{F}}}

\def\N{{\mathbb{N}}}

\def\R{{\mathbb{R}}}

\def\cB{{\mathcal B}}
\def\cC{{\mathcal C}}

\def\cE{{\mathcal E}}
\def\cF{{\mathcal F}}
\def\cH{{\mathcal H}}

\def\cM{{\mathcal M}}
\def\cP{{\mathcal P}}

\def\cS{{\mathcal S}}
\def\cU{{\mathcal U}}
\def\Ad{{\rm Ad \, }}
\def\Aut{{\rm Aut \, }}

\def\Ind{{\rm Ind }}

\newcommand{\norm}[1]{{\left\|{#1}\right\|}}

\newcommand{\abs}[1]{{\left|{#1}\right|}}
\newcommand{\scal}[1]{{\left\langle{#1}\right\rangle}}
\newcommand{\set}[1]{{\left\{{#1}\right\}}}

\newcommand{\dst}{\displaystyle}

\def\fig{ \centerline{Fig. \the\count200\global\advance\count200 by 1}}
\newtheorem{thm}{Theorem}[section]
\newtheorem{cor}[thm]{Corollary}
\newtheorem{lem}[thm]{Lemma}
\newtheorem{prop}[thm]{Proposition}

\newtheorem*{thm*}{Theorem}
\newtheorem*{lem*}{Lemma}
\newtheorem*{cor*}{Corollary}

\theoremstyle{definition}
\newtheorem{defn}[thm]{Definition}

\newtheorem{exs}[thm]{Examples}

\theoremstyle{remark}
\newtheorem{rem}[thm]{Remark}

\title[Comparison of norms of convolutors]{On the comparison of norms of   convolutors  associated
to noncommutative dynamics}

\author{Claire Anantharaman-Delaroche}
\address{D\'epartment de Math\'ematiques, Universit\'e d'Orl\'eans,
B. P. 6759, F-45067 Orl\'eans Cedex 2}
\email{claire.anantharaman@univ-orleans.fr}

\subjclass{Primary  46L55; Secondary 46L10, 22D25}

\keywords{Noncommutative dynamical systems, von Neumann algebras, standard form, amenability}

\begin{document}

\maketitle
\begin{abstract} To any action of a locally compact group $G$ on a pair $(A,B)$ of von Neumann algebras is canonically
associated a pair $(\pi_A^{\alpha}, \pi_B^{\alpha})$ of unitary representations of $G$. The purpose of this paper
is to provide results allowing to compare the norms of the operators $\pi_A^{\alpha}(\mu)$ and $\pi_B^{\alpha}(\mu)$
for bounded measures $\mu$ on $G$. We have a twofold aim. First to point out that several known facts in
ergodic and representation theory are indeed particular cases of general results about $(\pi_A^{\alpha}, \pi_B^{\alpha})$.
Second, under amenability assumptions, to obtain transference of inequalities that will  be useful
in noncommutative ergodic theory.
 
\end{abstract}

\section{Introduction}

Given a unitary representation $\pi$ of a locally compact group $G$ on a Hilbert space $\cH$ and
a bounded measure $\mu$ on $G$, we denote by $\pi(\mu)$
the operator $ \int_G \pi(s) d\mu(s)$ acting on $\cH$. It has long been known that  estimates of the spectral radius $r(\pi(\mu))$ or of
the norm $\norm{\pi(\mu)}$  give useful informations. 

A first  observation (often referred to as the ``Herz majorization principle''), asserts that for any locally compact group $G$, any closed
subgroup $H$ of $G$ and any {\em positive} bounded measure $\mu$ on $G$, one has 
$\norm{\lambda_G(\mu)}\leq \norm{\lambda_{G/H}(\mu)}$ (and therefore $r\big(\lambda_G(\mu)\big)\leq r\big(\lambda_{G/H}(\mu)\big)$ as well). 
Here $\lambda_G$ is the left regular representation and $\lambda_{G/H}$ denotes the  quasi-regular representation associated with $H$. More generally, given a representation $\pi$ of $H$, one has $\norm{\Ind_{H}^G\pi(\mu)} \leq
\norm{\lambda_{G/H}(\mu)}$, where $\Ind_{H}^G\pi$ is the representation induced of $\pi$, from $H$ to $G$ (see \cite{He}, \cite{EL}).

For a discrete group and a finitely supported symmetric probability measure $\mu$, 
the convolution operator $\lambda_G(\mu)$ had been investigated by Kesten \cite{Ke1} in connection with the random walk  defined by $\mu$.
 Kesten had already observed that $r\big(\lambda_G(\mu)\big)\leq r\big(\lambda_{G/H}(\mu)\big)$ when $H$ is a normal subgroup. Moreover in this seminal paper \cite{Ke1}, he proved that the normal subgroup $H$ is (what is now called) amenable if and only if there exists an adapted symmetric probability measure $\mu$ on $G$
({\it i.~e.} the support of $\mu$ generates the group $G$) such that $r\big(\lambda_{G/H}(\mu)\big) = r\big(\lambda_G(\mu)\big)$.

Note that in terms of operator algebras, these results concern the pair $$L^\infty(G/H)\subset L^\infty(G)$$ of
abelian von Neumann algebras, acted upon by left translations of $G$. 

In \cite{AD03}, to which this paper is a sequel, we considered $G$-actions on pairs $B\subset A$ of abelian von Neumann algebras. Our purpose in this paper is to deal
 more generally with $G$-actions on any pair $B\subset A$ of  von Neumann
algebras. In order to state the problems we are interested in, we need first to introduce some notations and definitions.
Let $A$ be a von Neumann algebra
and $\alpha$ a continuous homomorphism from a locally compact group $G$ into the group $\Aut\!(A)$ of
automorphisms of $A$. To such a {\it dynamical system} $(A, G, \alpha)$ is associated a unitary representation
$\pi^{\alpha}_A$ on the noncommutative $L^2$-space $L^2(A)$, well defined, up to
equivalence (see Section 2). Two particular examples are well known. The first one
is when $A = L^\infty(Y,m)$ is an abelian von Neumann algebra. In this case, $\alpha$
is an action on $Y$ which preserves the class of the measure $m$ and $\pi^{\alpha}_A$ is the corresponding
unitary representation in $L^2(Y,m)$. When $\alpha$ is the action of $G$ on $A=L^\infty(G/H)$ by 
translations, we get $\pi_A^{\alpha} = \lambda_{G/H}$. The second important example 
concerns the von Neumann algebra $A=\cB(\cH)$ of all bounded operators
on $\cH$ and, for $\alpha_s$, $s\in G$, the automorphism $T\mapsto \Ad\pi(s)(T) = \pi(s)T\pi(s)^*$ of $\cB(\cH)$,
where $\pi$ is a given   unitary representation  of $G$ on a Hilbert space $\cH$. Then the representation $\pi^{\alpha}_A$
is equivalent to the tensor product $\pi\otimes \overline{\pi}$ of $\pi$ with its conjugate $\overline{\pi}$.

 By a pair $(A,B)$ of von Neumann algebras, we mean that $B$ is a von Neumann
subalgebra of $A$. An {\it action $\alpha$ of $G$ on} $(A,B)$ is a dynamical system
$(A,G,\alpha)$ such that $B$ is globally $G$-invariant (we still denote by the same
letter the restricted action to $B$). Our first result is that the ``Herz majorization principle''  is valid
for every pair $(\pi^{\alpha}_A, \pi^{\alpha}_B)$. Namely:

\begin{thm*}{\rm (\ref{Herz0})} Let $\alpha$ be an action of $G$ on a pair $(A,B)$.  For every positive bounded measure $\mu$ on $G$ we have 
$$\norm{\lambda_G(\mu)}\leq \norm{\pi_A^{\alpha}(\mu)}\leq \norm{\pi_B^{\alpha}(\mu)}\leq \mu(G).$$
\end{thm*}

The proof is based on the fact that representations of the form $\pi^{\alpha}_A$ have enough $G$-positive
vectors, that is vectors $\xi$ such that $\scal{\xi,\pi^{\alpha}_A(s)\xi} \geq 0$ for all $s\in G$. Indeed, every normal
positive form \footnote{we shall denote by $A_*^{+}$ the cone of such forms} $\varphi$ on $A$ is represented by a well defined $G$-positive vector $\varphi^{1/2}$ in the space $L^2(A)$ of the
representation $\pi_A^{\alpha}$
and the set $L^2(A)^+$ of such vectors is a self-dual cone in $L^2(A)$. The second ingredient of the proof is the
inequality $\scal{\varphi^{1/2},\psi^{1/2}}_{L^2(A)}\leq\langle(\varphi_{ |_B})^{1/2},(\psi_{|_B})^{1/2}\rangle_{L^2(B)}$
for $\varphi,\psi\in A_{*}^+$ (where $\varphi _{|_B}, \psi_{|_B}$ are the restrictions of $\varphi$ and $\psi$ to $B$),
resulting from a variational formula due to Kosaki \cite{Ko}
(see Section 2 for notations and details).

In the case where $A = B\otimes M$ is the tensor product of two von Neumann algebras (with tensor product action)
we get $\norm{\big(\pi_{B}^\alpha\otimes \pi_{M}^\alpha\big)(\mu)}\leq \norm{\pi_{B}^\alpha(\mu)}$ for
every positive bounded measure $\mu$. In fact, this result remains true for any representation $\pi$
instead of $\pi_{M}^\alpha$ and any representation $\rho$ having a separating family of $G$-positive vectors
instead of $\pi_{B}^\alpha$ (see Theorem \ref{Herz2}). In particular, taking $\pi = \lambda_G$ and using
Fell's absorption principle, one gets $\norm{\lambda_G(\mu)} \leq \norm{\rho(\mu)}$. This inequality was proved by Pisier
 \cite{Pi} when $\rho$ is of the form $\pi\otimes \overline{\pi}$. Later,  Shalom observed in \cite{Sha} that this inequality
 holds for any representation having a non-zero $G$-positive vector.
 
Whereas the Herz majorization principle involves
 positivity properties, extensions of Kesten's result express amenability phenomena.
We say that the {\it action $\alpha$ of $G$ on $(A,B)$ is amenable} if there exists a norm one
projection\footnote{ a norm one projection $E$ from $A$ onto $B$ is also called 
a conditional expectation. It is automatically positive and satisfies $E(b_1 a b_2) = b_1 E(a) b_2$ for $a\in A$ and $b_1, b_2 \in B$} $E$ from $A$ onto $B$  such that $\alpha_s\circ E = E\circ \alpha_s$ for $s\in G$. Two special cases 
are particularly important. When $B = \C$, one says that the {\it $G$-action   on $A$ is  co-amenable} (or {\it amenable in the sense of Greenleaf}).  When
one considers the pair $(A\otimes L^\infty(G), A)$ with the tensor product action of the action on $A$ by left translations
on $L^\infty(G)$, one says that the {\it $G$-action on $A$ is amenable} (or {\it amenable in the sense
of Zimmer}).  Of course any action of an amenable group
on $(A,B)$ is amenable.

Let us consider an action $\alpha$ of $G$ on $(A,B)$. The existence of a 
{\em normal} $G$-equivariant conditional expectation from $A$ onto $B$ easily implies that $\pi^{\alpha}_B$ is a subrepresentation of $\pi^{\alpha}_A$.  It is therefore very natural to wonder whether
the amenability of the action implies that $\pi^{\alpha}_B$ is weakly contained in $\pi^{\alpha}_A$. We believe that this result
is true in general but we can only solve the problem in several particular cases.
Recall that a representation $\pi_1$ is said to be {\it weakly contained in a representation} $\pi_2$ (and we write $\pi_1\prec \pi_2$) if for every $f\in L^1(G)$ we have $\norm{\pi_1(f)}\leq \norm{\pi_2(f)}$ (or, equivalently, if $\norm{\pi_1(\mu)}\leq \norm{\pi_2(\mu)}$ for every bounded measure $\mu$ on $G$). 

Borrowing ideas used by A. Connes \cite{Co77} in order to show that injective von Neumann algebras are semi-discrete, we obtain:

 \begin{thm*}{\rm(\ref{faibleinc})} Let $\alpha$ be an amenable action of $G$ on a pair $(A,B)$ of von Neumann algebras. We assume that there is
 a faithful normal invariant state $\varphi$ on $B$. Then $\pi^{\alpha}_B$ is weakly contained in $\pi^{\alpha}_A$. 
 \end{thm*}

In presence of  ``enough'' {\em normal} conditional expectations from $A$
onto $B$, there is another approach aiming to approximate conditional expectations by normal ones. This gives  a stronger weak containment property. For instance we get:

\begin{thm*}{\rm (\ref{faibleinc1}, \ref{faibleinc2})}  Let $\alpha$ be an amenable action of $G$ on a pair $(A,B)$. We assume  either that $B$ is contained
in the centre $Z(A)$ of $A$ or that $A$  is a tensor product $B\otimes M$ of von Neumann algebras (with tensor product action).
 There exists a net $(V_i)$ of isometries
  from $L^2(B)$ into $L^2(A)$ such that for every $\xi\in L^2(B)$ one has
 $$\lim_i \norm{\pi_A^{\alpha}(s)V_i\xi - V_i \pi_B^{\alpha}(s)\xi} = 0$$
 uniformly on compact subsets of $G$. In particular, $\pi^{\alpha}_B$ is weakly contained in $\pi^{\alpha}_A$.
 \end{thm*}

There are more general statements (see Remark  \ref{partexp}). However we are mainly interested  in the case $(A\otimes L^\infty(G), A)$. As a consequence of
the previous theorem, we see that if an action $\alpha$ on $A$ is amenable, then 
for every bounded measure $\mu$  we have
\begin{equation}\label{regineg}
\norm{\pi_A^{\alpha}(\mu)}\leq \norm{\lambda_G(\mu)},
\end{equation}
these  inequality being  an equality when $\mu$ is positive.

Note that (\ref{regineg}) is a transference of norm estimates from the regular representation to $\pi_A^{\alpha}$,
a classical result when $A$ is abelian and $G$ amenable (see \cite{CW}). For an amenable action of $G$
on an abelian von Neumann algebra $A$, this inequality was obtained in \cite{Ku} when $G$ is discrete and
in \cite{AD03} for any locally compact group. It gives an upper bound for $\norm{\pi_A^{\alpha}(\mu)}$ only depending on $\mu$,
particularly useful in ergodic theory (for instance in the study of entropy, see \cite[Prop. 4.1]{Ne}).

As observed in \cite{AD03}, there is no hope for recovering in general the amenability of the action on
a pair $(A,B)$ from the weak containment $\pi_B^{\alpha}\prec \pi_A^{\alpha}$, even in the commutative setting. Let us consider for example
the case $A = L^\infty(G)$ and $B = L^\infty(G/H)$. We proved in \cite[Prop. 4.2.1, Cor. 4.4.5]{AD03} the equivalence of
the following three conditions:
\begin{itemize}
\item the action of $G$ on $(L^\infty(G),L^\infty(G/H))$ is amenable;
\item  $H$ is amenable;
\item $\lambda_{G/H} \prec \lambda_G$ and the trivial representation $\iota_H$ of $H$ is weakly
contained in the restriction of $\lambda_{G/H}$ to $H$.
\end{itemize}
For $H = SL(2,\R)$ and $G = SL(3,\R)$, one has $\lambda_{G/H} \prec \lambda_G$ although $H$ is not amenable
(see \cite[Section 4.2]{AD03}).

However, when $B = \C$, the situation is completely understood.
  Recall that  a bounded measure $\mu$ is said to
be {\it adapted} if the closed subgroup generated by the support ${\rm Supp}(\mu)$
of $\mu$ is $G$. Observe that for any representation $\pi$ of $G$ and any probability
measure $\mu$ on $G$, it is easily seen that $r(\pi(\mu)) = 1$  whenever the trivial representation $\iota_G$ of $G$
is weakly contained in $\pi$. However, the existence of an adapted probability measure $\mu$ on $G$ with $r(\pi(\mu)) = 1$ does not always imply $\iota_G \prec \pi$.

 \begin{thm*} {\rm (\ref{coamen})} Let $(A,G,\alpha)$ be a dynamical system. The following conditions are equivalent:
\begin{itemize}
\item[(i)] there exists a $G$-invariant state on $A$ (i.e. the action is co-amenable);
\item[(ii)] the trivial representation $\iota_G$ is weakly contained in $\pi^{\alpha}_A$;
\item[(iii)] there exists an adapted probability measure $\mu$ on $G$ with $r(\pi^{\alpha}_A(\mu)) = 1$.
\end{itemize}
\end{thm*}

The above theorem is well known when $A$ is an abelian von Neumann algebra. First, extending Kesten's and Day's results \cite{Ke1,Ke2,Day}, Derriennic and Guivarc'h \cite{DG} proved the equivalence between (ii) and (iii) when $A = L^\infty(G)$
 (see also \cite{BC}). In this case, the equivalence between (i) and (ii) is the Hulanicki-Reiter theorem (see \cite[Theorem 3.5.2]{Gre1}). Recall that $G$ is then said to be an amenable group.

 Next, the equivalence between (i) and (ii) was obtained by Eymard \cite{Ey} for $G$-homogeneous spaces $G/H$. Later, Guivarc'h \cite{Gui} proved that the previous theorem
  holds for any action on an abelian von Neumann algebra.
  
Another particular case of the above theorem concerns the von Neumann algebra $A=\cB(\cH)$ of all bounded operators
on $\cH$ and  $\alpha_s = \Ad \pi(s)$, $s\in G$, 
where $\pi$ is a given   unitary representation  of $G$ on  $\cH$. 
In this situation the equivalence of (i) and (ii) is due to Bekka \cite{Be} and the equivalence of the two last assertions
is a recent result of Bekka and Guivarc'h \cite{BG} \footnote{representations $\pi$ such that $\iota_G\prec \pi\otimes \overline{\pi}$ are called amenable}.

The equivalence between (i) and (ii) for any dynamical system $(A,G,\alpha)$ is proved by Kirchberg
in \cite[Sublemma 7.2.1]{Ki}. As an interesting consequence of this fact, note that the weak containment
of the trivial representation $\iota_G$ in $\pi^\alpha_{A}$ is independent of the topology of $G$.
 
This paper is organized as follows. We begin by recalling the basic facts to know on standard forms
of von Neumann algebras. In section 3  we prove noncommutative  Herz majorization theorems. In section 4
we study the weak containment property  $\pi^{\alpha}_B  \prec\pi^{\alpha}_A$
for an amenable action $\alpha$ of $G$ on $(A,B)$. Finally, in the last section we consider more specifically the cases
of amenable and coamenable actions.\\

For fundamentals of the theory of von Neumann algebras, we refer to \cite{Di2, KR}.   
{\em In the whole paper, we shall only consider  second countable locally compact groups, 
$\sigma$-finite measured spaces and von Neumann algebras with
separable preduals}, although these assumptions are not always necessary.

\section{Preliminaries on standard forms and noncommutative dynamical systems}

\subsection{Standard form of a von Neumann algebra}
Let $M$ be a von Neumann algebra. A {\it standard form of} $M$ is a normal faithful representation of $M$ in a Hilbert space
$\cH_M$ endowed with a conjugate linear isometric involution $J_M: \cH_M\to \cH_M$ and a self-dual cone $P_M\subset \cH_M$
such that
\begin{itemize}
\item $J_MMJ_M = M' $ (where $M'$ is the commutant of $M$ in $\cB(\cH_M)$;
\item $J_McJ_M = c^*$ for all $c\in M\cap M'$;
\item $J_M\xi = \xi$ for all $\xi\in P_M$;
\item $xJ_MxJ_M(P_M) \subset P_M$ for all $x\in M$.
\end{itemize}
Such a standard form $(M,\cH_M,J_M,P_M)$ exists and is unique, up to isomorphism. We refer to \cite{Ha75} for details
about this subject. Given  a faithful normal state $\varphi$ on $M$, one may take the standard representation to be the Gelfand-Naimark-Segal representation on $L^2(M,\varphi)$.
Denoting by $\xi_\varphi$ the unit of $M$, viewed in $L^2(M,\varphi)$, then $J_M$ is the antilinear isometry $J_\varphi$
given by the polar decomposition of the closure of $S_\varphi: x\xi_\varphi \mapsto x^*\xi_\varphi$, $x\in M$. Moreover,
$P_M$ is the norm closure $P_\varphi$ of $\{xJ_\varphi xJ_\varphi\xi_\varphi : x\in M\}$. 
 
 Usually we shall fix a standard form, denoted by
 $(M,L^2(M), J_M, P_M)$, or even $(M,L^2(M), J,P)$ for simplicity.
 The space $L^2(M)$ is ordered by the positive cone $P$.
This cone is self-dual in the sense that 
$$P = \{\xi\in L^2(M) : \scal{\xi,\eta} \geq 0, \forall \eta\in P\}.$$
 Recall also that
every element $\xi\in L^2(M)$ can be written in a unique way as $\xi = u\abs{\xi}$ where $\abs{\xi}\in P_M$ and $u$
is a partial isometry in $M$ such that $u^*u$ is the support of $\xi$. This decomposition is called the polar decomposition
of $\xi$. 

The Banach space $M_*$ of all normal forms on $M$ is the predual of $M$. A crucial fact is  that every normal positive form $\phi\in M_*^{+}$ can be uniquely written as
$\phi = \omega_\xi$, with $\xi\in P_M$ \footnote{where $\omega_\xi$ is the vector state $x\mapsto \scal{\xi,x\xi} = \omega_\xi(x)$}. It is also very suggestive
to denote $\phi^{1/2}$ this vector $\xi$.
 
We shall need a concrete description of $(L^2(M), J,P)$. We refer to \cite{Ha77} and \cite{Ter} for
the details concerning the following facts. We fix a concrete representation of $M$ on a Hilbert
space $\cH$. Let $\sigma : t\mapsto \sigma^{\psi}_t$ be the
modular automorphism group of a normal semi-finite faithful weight $\psi$ and let 
$$ M\rtimes_\sigma \R\subset \cB(L^2(\R)\otimes \cH)$$
 be the corresponding
crossed product. We denote by $\hat{\sigma}$ the dual action of $\R$ on $M\rtimes_\sigma \R$. Recall that $M\rtimes_\sigma \R$ has a canonical
normal semi-finite trace $\tau$ satisfying $\tau\circ\hat{\sigma}_t = e^{-t}\tau$ for all $t\in \R$. Following the point of
view of Haagerup, $L^p(M)$ is defined, for $p\geq 1$, as a  subspace of the $*$-algebra
$\cM(M\rtimes_\sigma \R)$ formed by the closed densely defined operators on $L^2(\R)\otimes \cH$,
affiliated with $M\rtimes_\sigma \R$, that are measurable with respect to $\tau$ (see \cite{Ter} or \cite{TaII}).
Namely, 
$$L^p(M) = \{ x\in\cM(M\rtimes_\sigma \R) : \hat{\sigma}_t(x) = e^{-t/p} x, \forall t\in \R\}.$$
In this picture, we have a description of $(M,L^2(M),J,P)$ as follows: $L^2(M)$ is defined, as just said, as a space of operators, its positive cone  is the cone of all positive operators in $L^2(M)$,
$J$ is the adjoint map, and the polar decomposition is the usual one. The spaces $L^1(M)$
and $L^\infty(M)$ are canonically isomorphic to the predual $M_*$ of $M$ and to $M$ respectively.

\begin{lem}\label{positiveL2} Let $(M, L^2(M), J, P)$ be a standard form of $M$. Then for $\xi,\eta \in L^2(M)$ we have
$$\abs{\scal{\xi,\eta}}^2 \leq \scal{\abs{\xi},\abs{\eta}}\scal{\abs{J\xi},\abs{J\eta}}.$$
\end{lem}

\begin{proof} The proof is exactly that of \cite[Lemma 2, p. 105]{Di2}. We use the above des\-cription of the standard form.
We need to introduce the linear functional
$$tr : h\in L^1(M) \mapsto \varphi_h(1),$$ where $\varphi_h$ denotes the normal linear form on $M$ associated with $h$.
We recall that if $h\in L^p(M)$ and $k\in L^q(M)$, with $1/p + 1/q = 1$, then
$tr(hk) = tr(kh)$. Moreover, for $\xi,\eta\in L^2(M)$ the scalar product is given by $\scal{\xi,\eta}= tr(\xi^*\eta)$.
The functional $tr$ plays the role of the usual trace on the space of trace-class operators on a Hilbert
space.

 Let $\xi = u\abs{\xi}$ and $\eta = v\abs{\eta}$ be the polar decompositions of $\xi$ and $\eta$ respectively. Then the polar decomposition of $J\xi$
is $\dst J\xi = u^*\big(u(JuJ)\abs{\xi}\big)$ since $J\nu = \nu$ for $\nu\in P$, so that $\abs{J\xi} = u\abs{\xi}u^*$.
Similarly, we have $\abs{J\eta} = v\abs{\eta}v^*$. Using the tracial property of $tr$ and the Cauchy-Schwarz inequality we 
proceed  as in  \cite[Lemma 2, p. 105]{Di2} to get 
$$\abs{\scal{\xi,\eta}} \leq tr(\abs{\xi}\abs{\eta})^{1/2} tr(\abs{J\eta}\abs{J\xi})^{1/2}.$$

\end{proof}

\subsection{Standard form of a pair of von Neumann algebras}
Now let $B$ be a von Neumann subalgebra of $A$ and let us examine some results relating the standard form $(B,L^2(B),J_B,P_B)$ of $B$ 
to that of $A$. 

Let us consider first a pair $(A,B)$ such that there exists a normal faithful conditional expectation $E$ from $A$ onto
$B$. Let us choose a faithful normal state $\psi$ on $B$ and set $\varphi = \psi\circ E$.
The Hilbert space $L^2(B,\psi)$ is canonically embedded into $L^2(A,\varphi)$ and the standard form of $B$
is obtained from that of $A$ by restriction to $L^2(B,\psi)$. Indeed, one checks
that $L^2(B,\psi)$ is stable under $J_\varphi$ and that $J_\psi$ is the restriction of $J_\varphi$ to $L^2(B,\psi)$.
Moreover, $P_\psi = P_\varphi \cap L^2(B,\psi)$ and the standard representation of $B$ into $L^2(B,\psi)$
is the restriction of the standard representation of $A$ into $L^2(A,\varphi)$ (see for instance \cite[page 130]{Str}).
For every $\xi,\eta \in L^2(B,\psi)$ and $a\in A$, one has $\scal{\xi, E(a)\eta} = \scal{\xi,a \eta}$. More generally, we shall
need the following lemma.

\begin{lem}\label{qu} Let $E$ be a normal conditional expection from $A$ onto $B$. There exists a unique positive isometry $q_E$
({\it i.e.} sending $L^2(B)^+$ into $L^2(A)^+$)
from $L^2(B)$ into $L^2(A)$ such that $\scal{\xi, E(a)\eta} = \scal{q_E(\xi), a q_E(\eta)}$ for all $a\in A$.
\end{lem}

\begin{proof}  The uniqueness of $q_E$ is immediate since for $\xi\in L^2(B)^+$, the normal form $\omega_\xi \circ E$ 
is uniquely implemented by $q_E(\xi)\in L^2(A)^+$. 
To prove the existence of $q_E$ we introduce the support $e$ of $E$, that is the smallest projection in $A$ with $E(1-e) = 0$. 
We have $e\in A \cap B'$ and $\theta : b\mapsto be$ is an isomorphism from $B$ onto $Be$. Through this isomorphism, we may view the standard form
of $B$ as represented into $L^2(Be)$ which is embedded in $L^2(eAe)$ by the previous remarks applied to the faithful normal
conditional expectation $a \mapsto eE(a)$ from $eAe$ onto $eB$. Now $L^2(eAe)$ is obviously included into $L^2(A)$. The composition of all these isometries
give the required $q_E : L^2(B) \to L^2(A)$.
\end{proof}

In the other direction we shall need the map $p$ from $L^2(A)^+$ into $L^2(B)^+$ defined by
$$\scal{\xi, b\xi} = \scal{p(\xi), bp(\xi)}$$
for all $b\in B$. In other terms, for $\varphi \in A_*^{+}$, we have $p(\varphi^{1/2}) = (\varphi_{|_B})^{1/2}.$

\begin{exs}\label{exap} (a) Assume that $B$ is a von Neumann subalgebra of $Z(A)$. We write $B = L^\infty(X, m)$ and
$A= \int_X^{\oplus} A(x) dm(x)$, so that $L^2(A) = \int_X^{\oplus} L^2(A(x)) dm(x)$ (see \cite{Su}).
Let $\xi = \int_X^{\oplus} \xi(x) dm(x)$ be an element of $L^2(A)^+$. Then $p(\xi)$ is the function $x\mapsto \norm{\xi(x)}_2$,
belonging to $L^2(X,m)^+$.

(b) Assume that $A$ is a finite von Neumann algebra, equipped with a faithful normal trace $\tau$ and let $E$ be the faithful normal conditional 
expectation from $A$ onto $B$ such that $\tau\circ E = \tau$. Then one has $A^+ \subset L^2(A,\tau)^+$ and 
$B^+ \subset L^2(B,\tau_{|_B})^+$.  For $a\in A^+$, one checks that $p(a) = E(a^2)^{1/2}$.

(c) Take $A = \cB(\cH)$ and $B= \C$. Let $\varphi = Tr(h\cdot)\in \cB(H)_*^{+}$ where $h$ is a positive trace-class operator on $\cH$.
Then $p(\varphi^{1/2}) =  p(h^{1/2}) = \norm{h^{1/2}}_2$.
 \end{exs}
 
 \begin{lem}\label{Kosaki} The map $p : L^2(A)^+ \to L^2(B)^+$ is a norm preserving homeomorphism and we have
 $$\scal{\xi,\eta} \leq \scal{p(\xi),p(\eta)}$$
 for all $\xi,\eta \in L^2(A)^+$.
 \end{lem}
 
 \begin{proof} Let $\varphi, \psi \in L^2(A)^+$ and let $\Delta_{\varphi,\psi}$ be the unique positive self-adjoint operator on $L^2(A)$
 such that $J\Delta_{\varphi,\psi}^{1/2} a \psi^{1/2} = a^*\varphi^{1/2}$ for all $a\in A$. Using 
 the formula
 $$\sqrt{\lambda} = \frac{1}{\pi}\int_0^{+\infty} \frac{\lambda}{\lambda + t}\frac{dt}{\sqrt{t}},$$
 we get
 \begin{align*}
 \scal{\varphi^{1/2},\psi^{1/2}} &= \scal{\Delta_{\varphi,\psi}^{1/2}  \psi^{1/2},\psi^{1/2}}\\
&= \frac{1}{\pi}\int_0^{+\infty} \scal{\Delta_{\varphi,\psi}\big(\Delta_{\varphi,\psi} +t)^{-1}\psi^{1/2} ,\psi^{1/2} }\frac{dt}{\sqrt{t}}.
\end{align*}
A quadratic interpolation method gives
 \begin{equation}\label{dimin}
 \scal{\Delta_{\varphi,\psi}\big(\Delta_{\varphi,\psi} +t)^{-1}\psi^{1/2} ,\psi^{1/2} } =
 \inf\set{\varphi(yy^*)/t + \psi(z^*z)}
 \end{equation}
where the infimum is taken on the pairs $(y,z)\in A^2$ such that $y+z = 1$ (see \cite[Lemma 2.1]{Ko} and the proof
of the Wigner-Yanase-Dyson-Lieb theorem \cite[Th. 5.2]{Ko}).
 We conclude by observing that the expression (\ref{dimin}) obviously increases when $\varphi$ and $\psi$
 are replaced by their restriction to $B$.

\end{proof}

\begin{rem} In examples \ref{exap} (a) and (c), the above lemma is immediately obtained by Cauchy-Schwarz
inequality. Applied to example \ref{exap} (b), this lemma gives the following inequality:
$$\forall x,y\in A^+,\quad \tau(xy) \leq \tau\big(E(x^2)^{1/2}E(y^2)^{1/2}\big).$$
\end{rem}

\subsection{Unitary implementation of automorphisms} Finally, let us recall a very important property of standard forms. Let $\Aut(M)$ be the automorphism group of the von Neumann algebra $M$ and let $(M,L^2(M), J,P)$ be a fixed standard form.
For every $\gamma\in \Aut(M)$ there is a unique  $u(\gamma)$ in the unitary group $\cU(L^2(M))$
such that $u(\gamma)(P) \subset P$, $Ju(\gamma) = u(\gamma) J$ and $\gamma(x) = u(\gamma) x u(\gamma)^*$ 
for all $x\in M$. This unitary is called
the {\it canonical implementation} of $\gamma$ (see \cite[Theorem 3.2]{Ha75}). 

The group $\Aut(M)$ acts on $M_*$ by $(\gamma,\varphi) \mapsto \varphi\circ \gamma^{-1}$. We equip it with the topology of pointwise norm convergence on $M_*$. Then the map $\gamma \mapsto u(\gamma)$ is a continuous homomorphism from
$\Aut(M)$ into the unitary group $\cU(L^2(M))$ equipped with the strong operator topology (\cite[Prop. 3.6]{Ha75}).

\begin{lem}\label{equiv} Let $(A,B)$ be a pair of von Neumann algebras and $\gamma\in \Aut(A)$ such that $\gamma(B) = B$.
We fix standard forms of $A$ and $B$ and denote by $u_A(\gamma)$ and $u_B(\gamma)$ the unitary implementations of 
$\gamma$ and of its restriction to $B$ respectively.
\begin{enumerate}
\item[(i)] We have $p\circ u_A(\gamma)(\xi) = u_B(\gamma)\circ p(\xi)$ for every $\xi\in L^2(A)^+$.
\item[(ii)] Let $E$ be a normal conditional expectation from $A$ onto $B$ and set $\gamma\cdot E = \gamma\circ E
\circ \gamma^{-1}$. Then $q_{\gamma\cdot E} = u_A(\gamma)\circ q_E \circ u_B(\gamma)^*$.
\end{enumerate}
\end{lem}

\begin{proof} Immediate. \end{proof}

\subsection{Representations defined by noncommutative dynamical systems} An {\it action of a locally compact group $G$ on a von Neumann algebra} $M$ is a continuous
homomorphism $s\mapsto \alpha_s$ from $G$ into $\Aut(M)$. We also say that
$(M,G,\alpha)$ is a {\it dynamical system}. We fix a standard form of $M$ and for $s\in G$,
we denote by $\pi^{\alpha}_M(s)$ the canonical unitary $u(\alpha_s)$ implementing $\alpha_s$.
 Then $\pi^{\alpha}_M$
is a unitary representation of $G$, that is a continuous homomorphism from $G$
into $\cU(L^2(M))$. Note that $\pi^{\alpha}_M$ is well defined, up to equivalence.

Observe that every $\xi \in P_M$ is {\it $G$-positive} in the sense that for all $s\in G$
we have $\scal{\xi, \pi^{\alpha}_M(s)\xi}\geq 0$. 

Note also that the set of  representations of the form $\pi^{\alpha}_M$ is stable
under direct sums and tensor products. Furthermore,  each representation $\pi^{\alpha}_M$ is equivalent to its conjugate. 

\begin{exs}\label{example} (a) Let $X$ be a standard Borel space with a Borel left $G$-action $(s,x)\in G\times X \mapsto sx\in X$.
When equipped with a $G$-quasi-invariant  measure $m$, we say that $(X,G,m)$ is a (non-singular)
{\it measured $G$-space}. To such a measured $G$-space is associated the dynamical system
$(L^\infty(X,m), G, \alpha)$ where $\alpha_s(f)(x) = f(s^{-1}x)$ for $f\in L^\infty(X,m)$, $s\in G$, $x\in X$.

We denote by $r$ (or $r_X$ in case of ambiguity) the {\it Radon-Nikod\' ym derivative} defined by
$$\forall s\in G, \forall f\in L^1(X,m), \quad \int_X f(s^{-1}x) r(x,s) dm(x) = \int_X f(x)dm(x).$$
Recall that $\big(L^\infty(X,m), L^2(X,m), J,L^2(X,m)^+\big)$  is a standard form of $L^\infty(X,m)$
(where $J$ is the complex conjugation and
$L^2(X,m)^+$ is the cone of non-negative functions in $L^2(X,m)$). The
unitary representation $\pi^{\alpha}_{L^\infty(X,m)}$ (rather denoted $\pi_X$) associated with
the dynamical system $(L^\infty(X,m), G, \alpha)$ is defined by
$$\pi_X(s)\xi(x) = \sqrt{r(x,s)} \xi(s^{-1}x)$$
for $\xi \in L^2(X,m)$ and $(s,x)\in G\times X$.

(b) Let $\pi$ be a representation of a locally compact group $G$ in a Hilbert space $\cH$.
We consider the dynamical system $(M,G,\alpha)$ where $M$ is the von Neumann algebra $ \cB(\cH)$
of all bounded operators on $\cH$ and $\alpha$ is the action such that $\alpha_s(T) = \pi(s) T \pi(s)^*$
for $T\in \cB(\cH)$ and $s\in G$. A standard form for $\cB(\cH)$ is $(\cB(\cH), \cH\otimes \overline{\cH}, J, (\cH\otimes \overline{\cH})^+)$
where $\cH\otimes \overline{\cH}$ is canonically identified with the Hilbert space of Hilbert-Schmidt operators,
$J$ is the adjoint operator (also described as $J: \xi\otimes\overline{\eta}\mapsto \eta\otimes\overline{\xi}$), and $(\cH\otimes \overline{\cH})^+$
is the cone of non-negative Hilbert-Schmidt operators. The canonical representation $\pi^{\alpha}_{\cB(\cH)}$
canonically associated with $(M,G,\alpha)$ is $\pi\otimes \overline{\pi}$. This representation  acts on the space
of Hilbert-Schmidt operators by
$$\pi\otimes \overline{\pi}(s)(T) = \pi(s)T\pi(s)^*.$$

(c) More generally, let $M$ be a von Neumann algebra and $s\mapsto \pi(s)$ be a continuous representation on $L^2(M)$
with $\pi(s) \in \cU(M)$ for all $s\in G$. Denote by $\alpha$ the corresponding action on $M$ by inner
automorphisms, that is $\alpha_s = \Ad \pi(s)$ for $s\in G$. Then $\pi_{M}^\alpha(s) =  \pi(s)J_M \pi(s)J_M$.

In the particular case where $M = L(G)$ is the group von Neumann algebra and $\pi = \lambda_G$ is the left
regular representation, $\pi_{M}^\alpha$ is the conjugation representation $\gamma_G$.
\end{exs}

\section{A non commutative ``Herz majorization principle''} 
 An {\it action of a locally compact group $G$ on the pair $(A,B)$ of von Neumann algebras} is a dynamical system 
$(A,G,\alpha)$ such that the von Neumann subalgebra $B$ is left globally invariant under the action. The restricted action of
$G$ on $B$ will still be denoted by $\alpha$.

\begin{thm}\label{Herz0} Let $\alpha$ be an action of the locally compact group $G$ on $(A,B)$. 
Then, for every probability measure $\mu$ on $G$, we have 
\begin{equation}\label{HHerz}
\norm{\pi_{A}^\alpha(\mu)}\leq
\norm{\pi_{B}^\alpha(\mu)}.
\end{equation}
\end{thm}

The proof uses the following well known way of computing norms.

\begin{lem}\label{estnorm} Let $T\in \cB(\cH)$ be a positive operator  on a Hilbert space $\cH$ and let $\cS\subset \cH$
be a separating family of norm one vectors ({\it i.e.} for $S\in \cB(\cH)$,  $S\xi = 0$ for all $\xi
\in \cS$ implies $S=0$). Then
$$\norm{T} = \lim_{n\to\infty}\sup_{\xi\in \cS}\omega_\xi(T^n)^{1/n}.$$
\end{lem}

\begin{proof} Let $\mu_\xi$ be the spectral measure of $T$ on the
spectrum $\sigma(T)\subset [0,\norm{T}]$, associated with $\xi$. Since the family $\cS$ is separating,
the union of the supports of $\mu_\xi$, $\xi\in \cS$, is dense into $\sigma(T)$. Given $\epsilon >0$, let
$\xi_0\in \cS$ be such that $\mu_{\xi_0}\big([\norm{T}-\epsilon,\norm{T}]\big) > 0$. We have
$$\norm{T} \geq \sup_{\xi\in \cS} \omega_\xi(T^n)^{1/n} \geq 
(\norm{T}-\epsilon)\mu_{\xi_0}\big([\norm{T}-\epsilon,\norm{T}]\big)^{1/n},$$
and the conclusion follows immediately.
\end{proof}

\begin{proof}[Proof of theorem \ref{Herz0}] We shall apply the previous lemma with $T= \pi_A^{\alpha}(\nu)$,
$\nu =\check{\mu}*\mu$, and for $\cS$ we take the set of norm one vectors in $L^2(A)^+$.
For every $n$, we have
$$\norm{\pi_{A}^\alpha(\mu)}^2 = \lim_{n\to\infty}\sup_{\xi\in \cS} \omega_\xi\big(\pi_A^{\alpha}(\nu^{*n})\big)^{1/n}$$
and
\begin{align*}
\omega_\xi\big(\pi_A^{\alpha}(\nu^{*n})\big) &= \int\scal{\xi,\pi_A^{\alpha}(t)\xi}d\nu^{* n}(t)\\
&\leq \int\scal{p(\xi),p\big(\pi_A^{\alpha}(t)\xi\big)}d\nu^{* n}(t)\\
&= \int\scal{p(\xi),\pi_B^{\alpha}(t)\big(p(\xi)\big)}d\nu^{* n}(t)\\
&=\omega_{p(\xi)}\big(\pi_B^{\alpha}(\nu^{*n})\big)\leq\norm{\pi_B^{\alpha}(\mu)}^{2n}
\end{align*}
by lemmas \ref{Kosaki} and \ref{equiv}.
The inequality (\ref{HHerz}) is then an immediate consequence of lemma \ref{estnorm}.
\end{proof}

 When $A$ is the tensor product of $B$ by another von Neumann algebra $C$, with a tensor product action
$\alpha = \beta\otimes \gamma$, it follows from theorem \ref{Herz0} that 
$$\norm{\pi_{B}^\beta\otimes\pi_C^{\gamma}(\mu)}
\leq \min\set{\norm{\pi_{B}^\beta(\mu)},\norm{\pi_{C}^\gamma(\mu)}}$$ for every probability measure $\mu$. 
In fact this is also a particular case of
the following more general result:

\begin{thm}\label{Herz2} Let $G$ be a locally compact group, $\rho$ a representation having a separating set $\cP$
of norm one $G$-positive vectors
and $\pi$ any representation. Then 
for every probability measure $\mu$ on $G$ we have
$$\norm{\big( \rho\otimes\pi\big)(\mu)} \leq \norm{\rho(\mu)}.$$
\end{thm}

\begin{proof} We  use again lemma \ref{estnorm} with $T = \big(\rho\otimes\pi\big)(\check{\mu}*\mu)$
and $$\cS = \set{\xi\otimes \eta : \xi \in \cP, \eta\in \cH(\pi), \norm{\eta} = 1}$$
 where $\cH(\pi)$
is the Hilbert space of the representation $\pi$.
  We have
$$\norm{\big( \rho\otimes\pi\big)(\mu)}^2=
\lim_{n\to\infty}\sup_{\xi\in \cS}\norm{\omega_\xi\big(( \rho\otimes\pi)(\nu^{*n})\big)}^{1/n}.$$
Observe that for $\xi\otimes\eta\in \cS$,
\begin{align*}
\omega_{\xi\otimes \eta}\big(( \rho\otimes\pi)(\nu^{*n})\big)&= 
\int \omega_{\xi}\big(\rho(t)\big)\omega_{\eta}\big(\pi(t)\big)
d\nu^{* n}(t)\\
&\leq\int \omega_{\xi}\big(\rho(t)\big)
d\nu^{* n}(t)\\
&= \omega_\xi\big(\rho(\nu^{*n})\big) \leq \norm{\rho(\mu)}^{2n},
\end{align*}
due to the positivity of $\omega_{\xi}\big(\rho(t)\big)$. 
The conclusion follows immediately.
\end{proof}

\begin{cor}{\rm (\cite[Lemma 2.3]{Sha})}\label{shalom} Let $\rho$ be a representation of $G$ having a non zero $G$-positive
vector. Then for every probability measure $\mu$ on $G$ we have
$$\norm{\lambda_G(\mu)}\leq \norm{\rho(\mu)}.$$
\end{cor}

\begin{proof} Let $\xi\in \cH(\rho)$ be a norm one $G$-positive vector for $\rho$. Then $\cP = \rho(G)\xi$
is a set of $G$-positive vectors. Let $K$ be the Hilbert subspace of $\cH(\rho)$ generated by $\cP$. It is $G$-invariant
and we denote by $\rho_{|_K}$ the representation of $G$ obtained by restriction. We apply theorem \ref{Herz2}
with $\rho_{|_K}$ instead of $\rho$ and $\pi = \lambda_G$. We have
$$\norm{(\rho_{|_K}\otimes \lambda_G)(\mu)}\leq \norm{\rho_{|_K}(\mu)} \leq\norm{\rho(\mu)}.$$
Moreover, a well known observation of  Fell \cite{Fe} says that the regular representation absorbs any other representation.
In particular $\rho_{|_K}\otimes \lambda_G$ is equivalent to a multiple of $\lambda_G$
and therefore $\norm{(\rho_{|_K}\otimes \lambda_G)(\mu)} = \norm{\lambda_G(\mu)}.$
\end{proof}

\begin{cor}{\rm (\cite{Le}, \cite{Pi})} Let $\pi$ be a representation of $G$. Then  for every probability measure $\mu$ on $G$ we have $\norm{\lambda_G(\mu)}\leq \norm{\big(\pi\otimes \overline{\pi}\big)(\mu)}$. Moreover, if
$\pi$ has a separating set of $G$-positive vectors, then 
$ \norm{\big(\pi\otimes\overline{\pi}\big)(\mu)}\leq \norm{\pi(\mu)}$.
\end{cor}

\begin{proof} The first inequality follows from Corollary \ref{shalom} and the second  from theorem \ref{Herz2}.
\end{proof}

\begin{rem} Let $U_1,\dots, U_n$ be $n$ unitary operators in a Hilbert space and let $c_1,\dots,c_n, d_1,\dots,d_n$
be $2n$ non negative real numbers. Let us denote by $g_1,\dots,g_n$ the generators of the free group $\F_n$.
As a particular case of the previous corollary one finds Pisier's inequality \cite{Pi}:
$$\norm{\sum_{i=1}^n \big(c_i\lambda_{\F_n}(g_i) + d_i\lambda_{\F_n}(g_i)^*\big)}
\leq \norm{\sum_{i=1}^n \big(c_i U_i\otimes\overline{U}_i + d_i U_{i}^*\otimes\overline{U}_{i}^*\big)}.$$
The left hand  side of this inequality is $2\sqrt{2n-1}$ when  $c_i=d_i=1$ for every $i$ (due to Kesten \cite{Ke1}).
It is equal to $2\sqrt{n-1}$ when $c_i = 1$ and $d_i=0$ for every $i$ (due to Akemann and Ostrand 
\cite{AO}).
\end{rem}

To conclude this section, let us explain how the inequality $\norm{\pi_A^{\alpha}(\mu)} \leq \norm{\pi_B^{\alpha}(\mu)}$, when $\mu$ is a probability measure and
$B\subset Z(A)$, is a particular case 
of the classical ``Herz majorization principle''. 

First we  need to  recall some definitions. Let $(X,m)$ be a measured
space and let $\cH = \{\cH(x): x\in X\}$
 be a $m$-measurable field of Hilbert spaces on $X$ (see \cite[Chap. II]{Di2}). We denote
by $L^2(\cH)= \int_{X}^{\oplus} \cH(x) dm(x)$ the direct integral Hilbert space. For $x,y \in X$,
the set of bounded linear maps from $\cH(x)$ to $\cH(y)$ will be denoted by $\cB\big(\cH(x),\cH(y)\big)$ and 
$\hbox{Iso}\big(\cH(x),\cH(y)\big)$ will be its subset of Hilbert space isomorphisms.  

\begin{defn}  Let $(X,G,m)$ be a measured $G$-space and let $\cH$ be as above. A {\it  (unitary) cocycle representation} of   $(X,G,m)$, acting on the measurable field
$\cH$, is a map 
$$U : (x,s)\in X\times G \mapsto U(x,s) \in\cB \big(\cH(s^{-1}x),\cH(x)\big)$$
 such that
\begin{itemize}
\item[(a)] for each $s \in G$, $U(x,s) \in \hbox{Iso}(\cH(s^{-1}x),\cH(x))$ for $m$-almost every $x$;
\item[(b)] for each $(s,t)\in G\times G$, $U(x,st) = U(x,s)U(s^{-1}x,t)$   for $m$-almost every
 $x\in X$;
\item[(b)] for every pair of measurable sections $\xi, \eta$ of $\cH$, and every $s\in G$
the map $x\mapsto \langle \eta(x),U(x,s)\xi(s^{-1}x)\rangle$ is measurable.
\end{itemize}
\end{defn}

To every cocycle representation $U$ of $(X,G,m)$ is associated a representation of $G$, called the {\it induced representation}
and denoted $\hbox{Ind} \,U$ (or $\hbox{Ind}_X U$ in case of ambiguity).
  Let us recall its definition. If $U$ acts on $\cH$, $\hbox{Ind} \,U$
 is the representation into
$L^2(\cH )$ defined by
$$(\hbox{Ind}\,U(s)\xi)(x) = \sqrt{r(x,s)} U(x,s)\xi(s^{-1}x)$$
for $\xi\in L^2(\cH)$ and $(x,s)\in X\times G$. This  extends the classical construction of
the representation  induced by a representation of a closed subgroup $H$, which amounts
to consider the left action of $G$ on $G/H$.

The fact that $\norm{\hbox{Ind} \,U(\mu)}\leq \norm{\pi_X(\mu)}$ when $\mu$ is a probability measure
is very easy to prove (see \cite[Prop 2.3.1]{AD03} for instance).

Now let $\alpha$ be an action of $G$ on a von Neumann algebra $A$, preserving a subalgebra $B = L^\infty(X,m)$
of $Z(A)$. By disintegrating $A$ with respect to $B$ we get $A = \int_X A(x) dm(x)$ and
$L^2(A) = \int_X L^2(A(x))dm(x)$ (see \cite{Su}).  We know by \cite[Th. 1]{Ma}
that the action of $G$ on $L^\infty(X,m)$ has a point realization. Therefore we may write $\pi_B^{\alpha}$ as
$$\pi_B^{\alpha}(s)\xi(x) = \sqrt{r(x,s)}\xi(s^{-1}x)$$
for $\xi\in L^2(X,m)$ and $(x,s)\in X\times G$. Thus we have $\pi_B^{\alpha} = \pi_X$. Moreover, by \cite[Prop. 1]{Gu} there is a cocycle representation $U_A : (x,s)\in X\times G \mapsto U_A(x,s)\in \cB\big(L^2(A(s^{-1}x)), L^2(A(x))\big)$ such that $\pi_A^{\alpha} = \hbox{Ind}\,U_A$.

\section{Representations associated with amenable pairs}

Let $(A,G,\alpha)$ be a dynamical system. It is easily checked that the representation $\pi_{A}^\alpha$
has a non-zero invariant vector if and only if there is a normal invariant state on $A$. More generally,
we have, in one direction:

\begin{prop} Let $\alpha$ be an action of $G$ on a pair $(A,B)$ of von Neumann algebras. Assume that there
exists a normal  $G$-equivariant conditional expectation $E$ from $A$ onto $B$. Then $\pi^{\alpha}_B$ is a subrepresentation of $\pi^{\alpha}_A$.
\end{prop}

\begin{proof}  This follows immediately from lemmas \ref{qu} and \ref{equiv}. Indeed, for an equivariant normal
conditional expectation $E$, the isometry $q_E$ of lemma \ref{qu} intertwines the representations 
$\pi_A^{\alpha}$ and $\pi_B^{\alpha}$.
\end{proof}

When the conditional expectation is not required to be normal, we are led to the following definition,
due to Zimmer \cite{Zi5} for pairs of abelian von Neumann algebras.

\begin{defn} We say that {\it an action of $G$ on a pair $(A,B)$ of von Neumann algebras is amenable} if there
exists an equivariant conditional expectation 
from $A$ onto $B$.
\end{defn}

Let us also recall the definitions of the two following important particular cases.

\begin{defn} Let $(A,G,\alpha)$ be a dynamical system.
\begin{enumerate}
\item[(i)] We say that the action is {\it co-amenable} if there is a $G$-invariant state on $A$.
\item[(ii)] We say that the action is {\it amenable} if there is a $G$-equivariant conditional expectation
from $A\otimes L^\infty(G)$ (with its usual tensor product action) onto $A$.
\end{enumerate}
\end{defn}

In particular, one of the definitions of amenability for a locally compact group $G$ is the co-amenability of the action on $L^\infty(G)$ by left translations.

 Let $\alpha$ be an action of $G$ on a pair $(A,B)$. Note that if the action on $B$ is amenable  then, by \cite[Prop. 2.5]{AD82}, the action on the pair $(A,B)$ is amenable whenever there exists a conditional expectation from
$A$ onto $B$.

In this section, we are interested in the following problem: 
\begin{itemize}
\item {\em does the amenability of the action $\alpha$ on $(A,B)$
imply that $\pi_B^{\alpha}$ is weakly contained in $\pi_A^{\alpha}$} ?
\end{itemize}
 We shall only give partial answers. First, we introduce some notations. For $f\in L^1(G)$,
  $a\in A$ and $\varphi \in A_*$ we set
 $$f*a = \int_G f(s)\, \alpha_s(a) ds,\quad (\varphi*f)(a) = \varphi(f* a).$$
 The left and right translated $s\!\cdot\! f$, $f\!\cdot\! s$ of $f$ are defined by
 $$(s\!\cdot\! f)(t) = f(s^{-1}t),\quad (f\!\cdot\! s)(t) = f(ts^{-1})\Delta(s)^{-1},$$
 where $\Delta$ is the modular function of $G$. Note that
 $f*\big(\alpha_s(a)\big) = (f\!\cdot\! s)*a$ and $\alpha_s(f*a) = (s\!\cdot\! f)*a$ for $a\in A$.
 
Let us now recall a useful equivalent definition of amenability using a  notion of invariant conditional expectation which is stronger than equivariance.

\begin{defn} Let us consider an action $\alpha$ of $G$ on a pair $(A,B)$ of von Neumann algebras. 
A {\it topologically invariant conditional expectation} is a conditional expectation $E: A\to B$ such that $E(f*a) = f*\big(E(a)\big)$
for every $f\in L^1(G)$ and every $a \in A$.
\end{defn}

 The following  result is well known
(see \cite{AD82} for instance):

\begin{prop} Let $G$ act on $(A,B)$. This action is amenable 
if and only if there exists a topologically invariant conditional expectation $E : A\rightarrow B$.
\end{prop}

\subsection{}Let us state a first answer to our problem.

 \begin{thm}\label{faibleinc} Let $\alpha$ be an amenable action of $G$ on a pair $(A,B)$ of von Neumann algebras. We assume that there is
 a faithful normal invariant state $\varphi$ on $B$. Then $\pi^{\alpha}_B$ is weakly contained in $\pi^{\alpha}_A$. 
 \end{thm}
 
 The proof uses the following key lemma inspired by the very simple proof given by Connes \cite{Co77} to show that an injective von Neumann algebra is semi-discrete.
 
 \begin{lem}\label{key} We keep the assumptions of the previous theorem. Given any compact subset $K$ of $G$
 and $\varepsilon >0$, there exists a normal state $\psi$ on $A$ such that
 $$\norm{\psi_{|_B} - \varphi} \leq \varepsilon, \quad \sup_{s\in K}\norm{\psi\circ\alpha_s - \psi}\leq\varepsilon.$$
 \end{lem}

 \begin{proof}   We follow the proof of \cite[Lemma 2]{Co77}.  We introduce the weakly compact convex set
 $$\cC = \set{x-\varphi(x)1 : x\in B, \norm{x}\leq \varepsilon^{-1}}.$$
 Let $E$ be a topologically invariant conditional expectation from $A$ onto $B$. The (non normal) state
 $\varphi\circ E$ belongs to polar set of the convex hull ${\rm co}(\cC\cup A^+)$, which is weakly closed. Therefore, using the bipolar
 theorem, we see that there is a net $(\psi_i)$ of normal states on $A$ such that $\lim_i \psi_i = \varphi\circ E$
 in the weak*-topology and $\norm{\psi_{i|_B}-\varphi}\leq \varepsilon$ for every $i$. 

 Now, we use the classical Day-Namioka convexity argument. We denote by $\cC'$ the convex
 set of normal states $\psi$ on $A$ such that $\norm{\psi_{|B}-\varphi}\leq \varepsilon$. Let $h_1,\cdots, h_k$ be fixed
 elements in $L^1(G)^+$ with $\int_G h_j(t)dt = 1$, $1\leq j\leq k$. For  $\psi\in \cC'$, we set
 $$b_j(\psi)=\psi * h_j - \psi.$$
 Let us denote by $\mathcal{C}''$ the range of $\cC'$ in the product
 $A_{*}^k$ by the map
 $$\psi\mapsto\big(b_1(\psi),\cdots,b_k(\psi)\big).$$
  Since $E$ is topologically invariant and $\varphi$ is invariant, we have
 $\big(\varphi\circ E\big)(h_j * a) =  \big(\varphi\circ E\big)(a)$ for $a\in A$. Therefore, we know that $(0,\cdots,0)$ belongs to the
 closure of $\mathcal{C}''$ in $A_{*}^k$ equipped with the product topology, 
 where we consider the weak topology on $A_*$. Since $\mathcal{C}''$
 is convex, we may replace this latter topology by the norm topology,
 using the Hahn-Banach separation theorem. It follows that there exists a net $(\psi_i)$
 in $\cC'$ such that for every $h\in L^1(G)^+$ with $\int_G h(t)dt = 1$ we have
  $$\lim_i\norm{\psi_i -\psi_i * h}_{A_*}=0 .$$
 
  Now we fix $h\in L^1(G)^+$
 such that $\int_G h(t)dt = 1$. Given $\eta >0$, we choose a neighbourhood $V$ of $e$
 in $G$ such that $\norm{h\!\cdot\!s -h}_1\leq \eta$ for $s\in V$. 
 Then we can find a finite number
 of elements $s_1,\dots, s_n$ in $G$ such that $K\subset \cup_{i=1}^n Vs_i $. We set $s_0 =e$
 and we choose $\psi\in \cC'$ satisfying 
$$ \norm{\psi -\psi*(h\!\cdot\! s_i)}_{A_*}\leq \eta$$
for  $0\leq i\leq n$.

Let $s \in K$ and choose $i$ such that $s\in Vs_i$. We have
\begin{align*}
 \norm{\psi*h -(\psi*h)\circ \alpha_s}_{A_*}
& \leq \norm{\psi*h -\psi}_{A_*}+
\norm{\psi - \psi*(h\!\cdot\! s)}_{A_*} \\
&\leq \eta + \norm{\psi - \psi*(h\!\cdot\! s_i)}_{A_*}+ \norm{\psi*(h\!\cdot\! s_i) - \psi*(h\!\cdot\! s)}_{A_*}\\
&\leq 2\eta + \norm{h - h\!\cdot\! (ss_{i}^{-1})}_{A_*} \leq 3\eta.
\end{align*}

To conclude, it suffices to  take $\eta= \varepsilon/3$  and to replace
$\psi$ by $\psi*h$.
 \end{proof}
 
 \begin{proof}[Proof of theorem \ref{faibleinc}] We fix $\varepsilon >0$ and a compact subset $K$ of $G$.
 Let $\psi$ be a normal state on $A$ as in lemma \ref{key}. We set $\xi_\varphi = \varphi^{1/2} \in L^2(B)^+$
 and $\xi_\psi = \psi^{1/2} \in L^2(A)^+$. Note that $\xi_\varphi = \pi_B^{\alpha}(s) \xi_\varphi$ for every $s\in G$
 since $\varphi$ is $G$-invariant. On the other hand, $\psi\circ \alpha_s$ corresponds to $\pi_A^{\alpha}(s) \xi_\psi$
 in $L^2(A)^+$. It follows from the Powers-St\o rmer inequality \cite[Lemma 2.10]{Ha75} that
 \begin{equation*}
 \norm{\xi_\psi - \pi_A^{\alpha}(s) \xi_\psi}_2^{2} \leq \norm{\psi - \psi\circ\alpha_s}_{A_*}.
 \end{equation*}
 
  Using the facts that $\norm{\xi_\psi - \pi_A^{\alpha}(s) \xi_\psi}_2\leq \sqrt{\varepsilon}$ and $\norm{\psi_{|_B} - \varphi}\leq
 \varepsilon$,  we get, for $s\in K$ and $b\in B$, 
  \begin{align*}
 \abs{\scal{b\xi_\psi, \pi_{A}^\alpha(s)b\xi_\psi}-  \scal{b\xi_\varphi, \pi_{B}^\alpha(s)b\xi_\varphi}}&= 
 \abs{\scal{\xi_\psi, b^*\alpha_s(b)\pi_{A}^\alpha(s)\xi_\psi} - \scal{\xi_\varphi, b^*\alpha_s(b)\xi_\varphi}}\\
 &\leq \abs{ \scal{\xi_\psi, b^*\alpha_s(b)\xi_\psi} - \scal{\xi_\varphi, b^*\alpha_s(b)\xi_\varphi}} + \norm{b}^2\sqrt{\varepsilon}\\
 &\leq \norm{b}^2(\varepsilon + \sqrt{\varepsilon}).
 \end{align*}

 Finally, we note that $B\xi_\varphi$ is dense in $L^2(B)$ since $\varphi$ is faithful. It follows from
 the above observations that every coefficient of the representation $\pi_B^{\alpha}$ is the limit, uniformly
 on compact subsets of $G$, of a net of coefficients of $\pi_A^{\alpha}$. Therefore, we have $\pi_B^{\alpha}\prec\pi_A^{\alpha}$.
 \end{proof}
 
 \subsection{} We now turn to situations where we can take advantage of the existence of sufficiently many normal conditional
 expectations.

 \begin{thm}\label{faibleinc1} Let $\alpha$ be an amenable action of $G$ on a pair $(A,B)$ of von Neumann algebras. We assume that $B$ is contained in the center $Z(A)$ of $A$. Then $\pi^{\alpha}_B$ is weakly contained in $\pi^{\alpha}_A$. More precisely, there exists a net $(V_i)$ of isometries
  from $L^2(B)$ into $L^2(A)$ such that for every $\xi\in L^2(B)$ one has
 $$\lim_i \norm{\pi_A^{\alpha}V_i\xi - V_i \pi_B^{\alpha}(s)\xi} = 0$$
 uniformly on compact subsets of $G$.
 \end{thm}

This last condition implies the weak containment property, as it is easily seen:

\begin{lem}\label{fortweak} Let $\pi$ (resp. $\rho$) be a representation on $\cH(\pi)$ (resp. $\cH(\rho)$). Assume 
the existence of a net $V_i$ of isometries  from $\cH(\rho)$ into $\cH(\pi)$ such that for every $\xi\in \cH(\rho)$
$$\lim_i \norm{\pi(s)V_i\xi - V_i \rho(s)\xi} = 0$$
 uniformly on compact subsets of $G$.
Then $\rho$ is weakly contained in $\pi$.
\end{lem}

\begin{proof}
 Let $\xi\in \cH(\rho)$ and  $f\in L^1(G)$. Since $V_i$ is an isometry, we have
 \begin{align*}
 \norm{\rho(f)\xi} &= \norm{V_i\rho(f)\xi}\\
 &\leq \norm{V_i\rho(f)\xi - \pi(f)V_i\xi}+ \norm{ \pi(f)V_i\xi}\\
& \leq \int_G \abs{f(s)}\norm{V_i\rho(s)\xi - \pi(s)V_i\xi}ds + \norm{\pi(f)}\norm{\xi}.
 \end{align*}
 Since $\lim_i\int_G \abs{f(s)}\norm{V_i\rho(s)\xi - \pi(s)V_i\xi}ds = 0$ it follows that  $\norm{\rho(f)\xi}\leq \norm{\pi(f)}\norm{\xi}$.  We conclude that $\rho$ is weakly contained in $\pi$.
 \end{proof}

 In order to prove  theorem \ref{faibleinc1}, we need some preliminaries.  We shall denote by $\cB_B(A,B)$ the Banach
 space of bounded maps $F$ from $A$ into $B$ that are $B$-linear in the sense that $F(b a ) = b F(a)$ for $a \in A$ and $b\in B$. This space is the dual of the quotient $A\widehat{\otimes}_B B_*$ of the projective tensor product
 $A\widehat{\otimes} B_*$ by the vector subspace generated by $\set{a\otimes b\varphi - ab \otimes \varphi : a\in A, b\in B, \varphi\in B_*}$. We denote by $\cB_B(A,B)_{1}^+$ the weak*-closed convex subset of positive elements
  $F\in \cB_B(A,B)$ with $F(1)\leq 1$. We want to introduce a weak*-dense convex subset $\cC$ in 
  $\cB_B(A,B)_{1}^+$, consisting of normal maps. To that purpose,
we disintegrate $A$ with respect to $B = L^\infty(X,m)$, that is
 we write $A = \int_{X}^\oplus A(x)dm(x)$ and $L^2(A) =  \int_{X}^\oplus L^2\big(A(x)\big)dm(x)$. Given a measurable
 section $\xi : x\mapsto \xi(x)\in L^2\big(A(x)\big)$ with $\norm{\xi(x)}_2 \leq 1$ almost everywhere, we denote by
 $\varpi_\xi$ the element of $\cB_B(A,B)_{1}^+$ such that $$\varpi_\xi(a)(x) = \scal{\xi(x), a(x)\xi(x)}_{L^2(A(x))}$$ for
 $a\in A$. Obviously, $\varpi_\xi$ is a normal element in $\cB_B(A,B)_{1}^+$ and we denote by $\cC$ the convex
 set of such maps $\varpi_\xi$. 
 
 \begin{lem}\label{dense} \begin{enumerate}
 \item[(i)] The set $\cC$ is weak*-dense in $\cB_B(A,B)_{1}^+$.
 \item[(ii)] Every conditional expectation from $A$ onto $B$ is the weak*-limit of a net of normal
 conditional expectations belonging to $\cC$.
 \end{enumerate}
 \end{lem}
 
 \begin{proof} Assume that there is an element $F\in \cB_B(A,B)_{1}^+$ that is not in the weak*-closure
 of $\cC$. Using the Hahn-Banach separation theorem, we find an element $\Phi = \sum_i a_i\otimes \varphi_i$
 in $A_{sa}\hat{\otimes}_B B_{*sa}$ (where the sum is finite) and $r\in \R$ with
 $$\scal{F, \Phi} > r\,\,\quad \hbox{and},\quad\,\, \forall \,\varpi_\xi\in \cC,\, \scal{\varpi_\xi, \Phi}\leq r.$$
 By polar decomposition, we may assume that $\varphi_i \geq 0$ for all $i$. Moreover, setting $\varphi = \sum_i
 \varphi_i$, thanks to the Radon-Nikod\'ym theorem we easily put $\Phi$ in the form $a\otimes\varphi$ with $a\in A_{sa}$.
 Let $e$ be a spectral projection of $a$ with $ae = a^+$, the positive part of $a$. Obviously,
 we have $\scal{F(a^+),\varphi} \geq \scal{F(a),\varphi}=\scal{F, \Phi} > r$. 
 
 On the other hand, observe that for every $\varpi_\xi \in \cC$
 the map $b\mapsto\varpi_\xi(eb)=\varpi_{\xi e}(b)$ still belongs to $\cC$. Therefore, writing $\varphi$ as $h\in L^1(X,m)^+$, we have
 $$\int_X h(x)\scal{\xi(x), a^+(x)\xi(x)}dm(x) = \scal{\varpi_{\xi e},\Phi}\leq r.$$
 
 Let $(\xi_n)$ be a sequence of measurable sections such that $\{\xi_n(x) : n\in \N\}$ is dense in the unit ball of $L^2(A(x))$ for almost every $x\in X$.   Let us fix $\varepsilon >0$.
 We may find a  measurable partition $(X_{n,\varepsilon})_n$ of $X$ such that for every $x\in X_{n,\varepsilon}$
 we have $\scal{\xi_n(x),a^+(x)\xi_n(x)} \geq \norm{a^+(x)}(1-\varepsilon)$. For $x\in X_{n,\varepsilon}$,
 we set $\xi(x)=\xi_n(x)$. Then we have
 $$r\geq \int_X h(x)\scal{\xi(x), a^+(x)\xi(x)}dm(x) \geq (1-\varepsilon)\int_X h(x)\norm{a^+(x)} dm(x).$$
 By letting $\varepsilon$ go to $0$ we get $r\geq \int_X h(x) \norm{a^+(x)} dm(x)$. 
 
 Since $F$ is positive, $B$-linear with $F(1) \leq 1$, we get $F(a^+)(x) \leq \norm{a^+(x)}$ {\it a.e.}
 and therefore 
 $$r< \scal{F(a^+),\varphi}=\int_X h(x) F(a^+)(x)dm(x) \leq \int_X h(x) \norm{a^+(x)} dm(x)\leq r.$$
 The contradiction thus obtained concludes the proof of (i).
 
 If $F(1) = 1$, it is easy  to see that we may approximate $F$ by elements $\varpi_{\xi}$ with
 $\norm{\xi(x)}_2 = 1$ {\it a.e.}
 \end{proof}

 We shall denote by $\cE(A,B)\subset \cB_B(A,B)_{1}^+$
the subset of all normal conditional expectations from $A$ onto $B$.

 \begin{lem}\label{fstep}  Let $\alpha$ be an action of $G$ on $(A,B)$ where $B\subset Z(A)$.
The following conditions are equivalent:
 \begin{itemize}
 \item[(i)] There exists a  $G$-equivariant conditional expectation from $A$ onto $B$.
 \item[(ii)] There exists a net $(\Phi_i)$ of normal conditional expectations from $A$ onto $B$ such that for $\varphi\in B_*$, $f\in L^1(G)$ and $a\in A$ we have
 $$\lim_i \scal{\varphi,f*\big(\Phi_i(a)\big) -\Phi_i(f*a)} = 0.$$
 \end{itemize}
\end{lem}

\begin{proof} (i) $\Rightarrow$ (ii). Let $E$ be a topologically invariant
conditional expectation and let $(\Phi_i)$ be a net in $\cE(A,B)$ such that $\lim_i \Phi_i = E$ in the weak*-topology,
whose existence was proved in lemma \ref{dense}.
Assertion (ii) follows immediately from the invariance
 of $E$.
 
The converse is also obvious.
\end{proof}

 To go further, let us introduce some more notations. For $s\in G$ and $F\in \cB_B(A,B)$ we  set $s\cdot F = \alpha_s\circ F\circ \alpha_{s^{-1}}$. Note that
 $s\cdot F \in \cE(A,B)$ whenever $F\in \cE(A,B)$. Finally, for $F\in \cE(A,B)$ and $f\in L^1(G)$,
 we define $f*F\in \cE(A,B)$ by
 $$\forall a\in A, \quad (f*F)(a) = \int_G f(s) (s\!\cdot\! F)(a) ds.$$
 
\begin{lem}\label{sstep} Let $\alpha$ be an action of $G$ on a pair $(A,B)$. The following conditions are equivalent:
\begin{itemize}
 \item[(i)] There exists a net $(\Phi_i)$ in $\cE(A,B)$ such that for every $\varphi\in B_*$,  $f\in L^1(G)$ and $a\in A$
 we have 
 $$\lim_i\scal{\varphi, f*\big(\Phi_i(a)\big) - \Phi_i(f*a)} = 0.$$
 \item[(ii)] There exists a net $(\Phi_i)$ in $\cE(A,B)$ such that for every $\varphi\in B_*$ and every $f\in L^1(G)$
 with $\int_G f(s) ds = 1$ we have 
 $$\lim_i\norm{\varphi\circ\big(f*\Phi_i - \Phi_i\big)}_{B_*} = 0.$$
 \item[(iii)] For every compact subset $K$ of $G$, every finite subset $\cF$ of $L^2(B)$ and every
 $\varepsilon >0$, there exists $\Phi\in \cE(A,B)$ such that
 $$\sup_{(s,\xi)\in K\times \cF} \norm{\pi_A^{\alpha}(s)q_{\Phi}\xi -q_{\Phi}\pi_B^{\alpha}(s)\xi}_{2} \leq \varepsilon.$$
 \item[(iv)] There exists a net $(\Phi_i)$ in $\cE(A,B)$ such that for every $f\in L^1(G)$ and $\xi \in L^2(B)$ we have $$\lim_i\int f(s) \norm{\pi_A^{\alpha}(s)q_{\Phi_i}\xi -q_{\Phi_i}\pi_B^{\alpha}(s)\xi}_{2} ds = 0.$$
 \end{itemize}
 \end{lem}

 \begin{proof} (i) $\Rightarrow$ (ii) Let $(\Phi_i)$ as in the statement of (i). We have
 \begin{align*}
\scal{\varphi, f*\big(\Phi_i(a)\big) - \Phi_i(f*a)} 
 &= \int_G f(s)\scal{\varphi, \alpha_s\circ\Phi_i(a) - \Phi_i\circ\alpha_s(a)}ds\\
 &= \int_G f(s)\scal{\varphi\circ\alpha_s, \Phi_i(a) - \alpha_{s^{-1}}\circ\Phi_i\circ\alpha_s(a)}ds.
 \end{align*}
 
Note that $s\mapsto f(s)\varphi\circ\alpha_s$ is in $L^1(G,B_*)$ and that elements of this form
generate $L^1(G,B_*)$. It follows that 
 $$\lim_i \int_G \scal{h(s),\Phi_i(a) - \alpha_{s^{-1}}\circ\Phi_i\circ\alpha_s(a)}ds =0$$
 for every  $h\in L^1(G,B_*)$ and $a\in A$ .

 Now, we use again the  Day-Namioka convexity argument. Let $h_1,\cdots, h_k$ be fixed
 elements in $L^1(G,B_*)$. For $1\leq j\leq k$ and $\Phi\in \cE(A,B)$, we set
 $$b_j(\Phi)= \int_G h_j(s)\circ\big(\Phi-s\!\cdot\!\Phi\big)ds \in A_*.$$
 Let us denote by $\mathcal{C}'$ the range of $\cE(A,B)$ in the product
 $A_{*}^k$ by the map
 $\Phi\mapsto\big(b_1(\Phi),\cdots,b_k(\Phi)\big)$. We know that $(0,\cdots,0)$ belongs to the
 closure of $\mathcal{C}'$ in $A_{*}^k$ equipped with the product topology, 
 where we consider the weak topology on $A_*$. Since $\mathcal{C}'$
 is convex, we may replace this latter topology by the norm topology. Therefore there exists a net $(\Phi_i)$
 in $\cE(A,B)$ such that for every $f\in L^1(G)$ with $\int_G f(s)ds = 1$ and every
 $\varphi\in B_*$, we have
 
  $$\lim_i\norm{\varphi\circ\big(\Phi_i -f*\Phi_i}_{A_*}\\
 =\lim_i  \norm{\int_G f(s)\varphi\circ\big( \Phi_i -s\!\cdot\!\Phi_i\big)ds}_{A_*} =0 .$$
 
 (ii) $\Rightarrow$ (iii). Let $K$ be a compact subset of $G$ and $\mathcal{F}$ a finite subset of $L^2(B)^+$. Let $f\in L^1(G)^+$
 such that $\int_G f(s)ds = 1$. We argue as in the proof of lemma \ref{key} to show
 that, given $\eta >0$, there exists $\Psi\in \cE(A,B)$ such that 
$$\sup_{(s,\xi)\in K\times \cF}
 \norm{\omega_\xi\circ\Big(s\!\cdot\!(f*\Psi)-f*\Psi\Big)}_{A_*} 
\leq \eta.$$

We take $\Phi = f*\Psi$. Using the Powers-St\o rmer inequality 
and lemma \ref{equiv} we get, for $s\in K$ and $\xi\in \cF$,
\begin{align*}
\norm{\pi^{\alpha}_A(s^{-1})q_{\Phi}\xi - q_{\Phi}\pi_B^{\alpha}(s^{-1})\xi}^{2}_2 
&= \norm{q_{\Phi}\xi - \pi^{\alpha}_A(s)q_{\Phi}\pi_B^{\alpha}(s^{-1})\xi}^{2}_2\\
&\leq \norm{\omega_\xi\circ\big(\Phi - s\!\cdot\! \Phi\big)}_{A_*} .
\end{align*}

To conclude, it suffices to replace $K$ by $K^{-1}$ and to take $\eta=\varepsilon^2$.

 (iii) $\Rightarrow$ (iv) is obvious. 
 It is not difficult to show that (iv) implies (i) and we skip the
 proof.
 \end{proof}

\begin{rem} \label{partexp} An inspection of the above proof shows that when $\alpha$ is an amenable action on $(A,B)$ with
$B\subset Z(A)$, one may take the $\Phi_i$'s in the convex set $\cC$ of lemma \ref{dense}. Using the cocycle
representation $U_A$ introduced at the end of Section 3, and taking $\Phi_i = \varpi_{\xi_i}$ we get,
for $\eta \in L^2(B) = L^2(X,m)$,
\begin{align*}
\norm{q_{\Phi_i}\pi_B^{\alpha}(s^{-1})\eta -\pi_A^{\alpha}(s^{-1})q_{\Phi_i}\eta}_{2}^2 
&= \norm{\pi_A^{\alpha}(s)q_{\Phi_i}\pi_B^{\alpha}(s^{-1})\eta -q_{\Phi_i}\eta}_{2}^2 \\
&=\int_X \norm{
U_A(x,s)\xi_i(s^{-1}x)\eta(x) -\xi_i(x)\eta(x)}_2^{2} dm(x)\\
& =\int_X |\eta(x)|^2\norm{U_A(x,s)\xi_i(s^{-1}x) - \xi_i(x)}_2^{2}dm(x).
\end{align*}
 Now by lemma \ref{sstep} (iv), we see that $\alpha$ is amenable if and only
if there exists a sequence $(\xi_n)$ of sections of the Hilbert bundle $(L^2\big(A(x)\big)_{x\in X}$ with
$\norm{\xi_n(x)}_2 = 1$ almost everywhere, such that
$$\lim_n \int_{X\times G} f(x,s)\norm{U_A(x,s)\xi_n(s^{-1}x) - \xi_n(x)}_2^{2} dm(x) ds = 0$$
for every $f\in L^1(X\times G)$. Expressed in term of groupoid, this is equivalent to the fact that
the representation $U_A$ of the measured groupoid $X\rtimes G$ weakly contains the trivial representation
(see \cite{AD04} for details).
\end{rem}

  \begin{proof}[Proof of theorem \ref{faibleinc1}] Immediate consequence of lemmas \ref{fstep} and \ref{sstep}.
  \end{proof}


 
\begin{thm}\label{faibleinc2} Let $\alpha$ be a tensor product action on $A = B\otimes M$. Assume that there is an equivariant
conditional expectation from $A$ onto $B$. Then the conclusions of theorem \ref{faibleinc1} hold. In
particular for every probability measure $\mu$ on $G$ we have
$\norm{\pi_B^{\alpha}(\mu)} = \norm{\big(\pi_B^{\alpha}\otimes \pi_M^{\alpha}\big)(\mu)}$.
\end{thm}

\begin{proof} We first observe that by restriction there is an equivariant conditional expectation
from $Z(B)\otimes M$ onto $Z(B)$. We write $Z(B)$ as $L^\infty(X,m)$ and we disintegrate the representation
$\pi_B^{\alpha}$ so that for $\xi =\int_X^{\oplus} \xi(x) dm(x) \in \int_X^{\oplus} L^2(B(x)) dm(x)$ we have (see the end of Section 3),
$$\pi_{B}^{\alpha}(s)\xi(x) =  \sqrt{r(x,s)}U_B(x,s)\xi(s^{-1}x).$$

Using Theorem \ref{faibleinc1} and remark \ref{partexp}, we get a net $(\Phi_i)$ of conditional expectations
from $Z(B)\otimes M$ onto $Z(B)$, of the form $\Phi_i = \varpi_{\xi_i}$ (where $\xi_i : X \to L^2(M)$ is a measurable map with $\norm{\xi_i(x)}_2 = 1$ almost everywhere), such that for every  $\eta \in L^2(Z(B))^+$
 we have 
 $$\lim_i \norm{\pi_{Z(B)\otimes M}^{\alpha}(s)q_{\Phi_i}\eta -q_{\Phi_i}\pi_{Z(B)}^{\alpha}(s)\eta}_{2}  = 0$$
 uniformly on compact subsets of $G$.
 
 Each linear isometry $q_{\Phi_i} : L^2(Z(B))\to L^2(Z(B))\otimes L^2(M)$ extends to an
 isometry $q_{\Phi_i} : L^2(B)\to L^2(B)\otimes L^2(M)$ by setting
 $$q_{\Phi_i}\eta(x) = \eta(x)\otimes \xi_i(x)$$
   for $\eta= \int_X^{\oplus} \eta(x) dm(x) \in L^2(B)$.
   
     A straightforward computation shows that for $\eta \in L^2(B)$,
 $$\norm{\pi_{B\otimes M}^\alpha(s)q_{\Phi_i}\eta -q_{\Phi_i}\pi_{B}^\alpha(s)\eta}_{2}
 =\norm{\pi_{Z(B)\otimes M}^\alpha(s)q_{\Phi_i}\abs{\eta} -q_{\Phi_i}\pi_{Z(B)}^\alpha(s)\abs{\eta}}_{2}$$
 where we denote by $\abs{\eta}$ the element $x\mapsto \norm{\eta(x)}_2$ of $L^2(Z(B))$.
 This ends the proof.
 \end{proof}
 
 \begin{rem} We may more generally use the same kind of techniques for any action $\alpha$
 on $B\otimes M$ leaving $B\otimes 1$ invariant even if the action is not a tensor product action. We may
 even deal  with an action on a pair $(B\otimes M, B\otimes N)$ where $N\subset Z(M)$, such
 that $B\otimes N$ is globally invariant under $\alpha$. As a consequence, using the structure theory
 of type $I$ von Neumann algebras, on gets the following result:
 
 \begin{thm} Let $\alpha$ be an amenable action of $G$ on a pair $(A,B)$ where $B$ is a type $I$
 von Neumann algebra such that $Z(B)\subset Z(A)$. Then there exists a net   $(V_i)$ of isometries
  from $L^2(B)$ into $L^2(A)$ such that for every $\xi\in L^2(B)$ one has
 $$\lim_i \norm{\pi_A^{\alpha}V_i\xi - V_i \pi_B^{\alpha}(s)\xi} = 0$$
 uniformly on compact subsets of $G$. In particular  $\pi^{\alpha}_B$ is weakly contained in $\pi^{\alpha}_A$.
\end{thm}

Since we are mainly interested in amenable and co-amenable actions, we shall not give the rather
tedious proof. In fact one would be more interested in deciding whether the above theorem is true
when $\alpha$ is an amenable action on $(A,B)$, under the assumption that there are ``enough'' normal
conditional expectations from $A$ onto $B$.
\end{rem}

Recall that a group $G$ is said to have {\it property $T$} if every of its representations that weakly
contains the trivial representation $\iota_G$ actually contains $\iota_G$ as a subrepresentation, that is has a non-zero $G$-invariant vector. It follows that for such groups, every dynamical system having an invariant state, {\it i.e.} $\iota_G\prec
\pi_A^{\alpha}$, has a normal invariant state, {\it i.e.} $\iota_G\leq\pi_A^{\alpha}$. More generally, we have:

\begin{thm} Let $G$ be a locally compact group having property $T$ and let $\alpha$
be an amenable action on a pair $(A,B)$ of von Neumann algebras. We assume that $B$ is contained in the 
centre of $A$ and that the action on $B$ is ergodic and leaves invariant a normal faithful state. Then there exists an equivariant
normal conditional expectation from $A$ onto $B$. In particular, $\pi_B^{\alpha}$ is a subrepresentation
of $\pi_A^{\alpha}$.
\end{thm}

\begin{proof} We write $B = L^\infty(X,m)$ where $m$ is an invariant probability measure. Since $G$ 
has property $T$ and preserves the finite measure $m$, one knows that the measured groupoid $X\rtimes G$
has property $T$ (see \cite[Cor. 5.16]{AD04} for instance). By remark \ref{partexp}, its  representation $U_A$ weakly
contains the trivial one. By definition of property $T$ for a measured groupoid, the trivial representation is actually contained in $U_A$. This means that there exists
a section $\xi : x \mapsto \xi(x) \in L^2\big(A(x)\big)$ with $\norm{\xi(x)}_2 = 1$ almost everywhere,
such that $U_A(x,s)\xi(s^{-1}x) = \xi(x)$ almost everywhere on $X\times G$. Then $E = \varpi_{\xi}$
is  a normal equivariant conditional expectation from $A$ onto $B$.
\end{proof}

\subsection{} We now mention another  positive answer to the problem considered in this section.
Let $B$ be a von Neumann algebra, $G$ a locally compact group and $\alpha$ the $G$-action on $B$ associated to a representation of $\pi : G \to  \cU(B)$ (see example \ref{example} (c)).
This action extends  to the action $\alpha : s\mapsto \Ad \pi(s)$ on $A= \cB\big(L^2(B)\big)$.
Observe that the amenability of the $G$-action on $(A,B)$ is equivalent to the injectivity of $B$,
that is to the existence of a norm one projection from $\cB\big(L^2(B)\big)$ onto $B$.

\begin{prop}\label{UH} Let $\alpha$ be the $G$-action  on $(A,B)$ defined above. Assume that $B$
is injective. Then we have $\pi_{B}^\alpha \prec \pi_{A}^\alpha$, that is $\pi J_B\pi J_B \prec \pi\otimes \overline{\pi}$.
\end{prop}

\begin{proof} By  the result asserting that an injective von Neumann algebra is semi-discrete (see \cite{Co, Co77})
we have, for every $a_i, b_i \in B$, $i=1,\dots,n$, 
$$\norm{\sum_{i=1}^n a_iJ_B b_i J_B}\leq \norm{\sum_{i=1}^{n} a_i\otimes \overline{b_i}}.$$
In particular, for $f_i,g_i$ in $L^1(G)$ we get
$$\norm{\sum_{i=1}^n \pi(f_i) J_B\pi(g_i)J_B}\leq \norm{\sum_{i=1}^{n} \pi(f_i)\otimes \overline{\pi}(g_i)}.$$
It follows that
\begin{equation}\label{equGG}
\norm{\int_{G\times G} h(s,t) \pi(s)J_B \pi(t)J_B ds dt} \leq \norm{\int_{G\times G} h(s,t) \pi(s)\otimes \overline{\pi}(t) ds dt}
\end{equation}
for $h\in L^1(G\times G)$.
Let $\mu$ be a bounded measure on $G$ and denote by $\nu$ its image by the diagonal map $s\mapsto (s,s)$.
Moreover let us consider an approximate unit $(\varphi_i)$ of $L^1(G)$. Applying the inequality (\ref{equGG})
to $h = (\varphi_j \otimes \varphi_j)*\nu$ we get
$$\norm{\pi(\varphi_j) J_B\pi(\varphi_j)J_B\int_G \pi(s)J_B\pi(s) d\mu(s)}\leq \norm{\pi(\varphi_j)\otimes\overline{\pi}(\varphi_j)\int_G \pi(s)\otimes\overline{\pi}(s)d\mu(s)}$$
from which we easily get
$\norm{(\pi J_B\pi J_B)(\mu)}\leq \norm{( \pi\otimes \overline{\pi})(\mu)}$.
\end{proof}

\begin{rem} If we apply this proposition to  the group von Neumann algebra $B= L(G)$ and to the 
representation $\pi = \lambda_G$, we get that whenever $L(G)$ is injective (for instance when $G$ is almost connected by \cite[Cor. 6.7]{Co}),
then  the conjugation representation $\gamma_G$ of $G$ is weakly contained in $\lambda_G \otimes \overline{\lambda}_G$ and therefore in $\lambda_G$.
Note that when the reduced $C^*$-algebra of $G$ is nuclear, it has been proved by Kaniuth \cite{Ka} that $\gamma_G$ is weakly contained in the direct
sum of the representations $\pi\otimes \overline{\pi}$ where $\pi$ ranges over the reduced dual of $G$.
\end{rem}

\section{Amenable and coamenable actions}

\begin{prop}\label{weakcont2} Let $(A,G,\alpha)$ be an amenable dynamical system. We have
$$  \pi^{\alpha}_{Z(A)} \prec \pi^{\alpha}_A \prec \lambda_G.$$
 In particular for every probability measure $\mu$
on $G$ we have
$$\norm{\lambda_G(\mu)} = \norm{\pi^{\alpha}_A(\mu)} = \norm{\pi^{\alpha}_{Z(A)}(\mu)}.$$
\end{prop}

\begin{proof} By \cite[Cor. 3.6]{AD82}, we know that the action of $G$ on $Z(A)$ is amenable. 
Moreover, it follows from \cite[Prop. 2.5]{AD82} that the action of $G$ on the pair $(A,Z(A))$
is amenable. By theorems \ref{faibleinc1} and \ref{faibleinc2}, we have respectively
$\pi^{\alpha}_{Z(A)} \prec \pi^{\alpha}_A$ and $\pi^{\alpha}_A\prec \pi^{\alpha}_A\otimes \lambda_G$.
Then the  conclusion follows from Fell's absorption principle.

The second part of the proposition follows from theorem \ref{Herz0}.
\end{proof}

\begin{rem} The property $ \pi^{\alpha}_A \prec \lambda_G$ means that for every bounded  measure $\mu$
on $G$ we have
$\norm{\pi^{\alpha}_A(\mu)}\leq \norm{\lambda_G(\mu)}$. It is a transference property of norm estimates, in the style
of the ones that prove to be so useful in classical harmonic analysis and ergodic theory (see \cite{CW}). In the noncommutative
setting, one can also establish $L^p$-transference inequalities and apply them to prove ergodic theorems.
This is the subject of a forthcoming paper.
\end{rem}

\begin{thm}\label{coamen} Let $(A,G,\alpha)$ be a dynamical system. The following conditions are equivalent:
\begin{itemize}
\item[(i)] there exists a $G$-invariant state on $A$ ({\it i.e.} the action is coamenable);
\item[(ii)] the trivial representation $\iota_G$ is weakly contained in $\pi^{\alpha}_A$;
\item[(iii)] there exists an adapted probability measure $\mu$ on $G$ with $r(\pi^{\alpha}_A(\mu)) = 1$.
\end{itemize}
\end{thm}

\begin{proof}
 (i) $\Rightarrow$ (ii) is a particular case of theorem \ref{faibleinc1} where we take $B=\C$
 (in this case, all technical difficulties disappear and the proof is indeed straightforward
by usual convexity arguments).

 (ii) $\Rightarrow$ (iii) is obvious. In fact, if (ii) holds, one easily sees that $1$ is an approximate eigenvalue
of $\pi^{\alpha}_A(\mu)$ for any probability measure $\mu$ on $G$.

To show (iii) $\Rightarrow$ (i) we follow the lines of the proof of Theorem 1 in \cite{BG} (or
of \cite[Th\'eor\`eme]{DG}), that we reproduce for the reader's convenience. First, since $\pi^{\alpha}_A(\mu)$ is a contraction of spectral
radius $1$, there exist a complex number $c$ with $\abs{c} = 1$ and a sequence
$(\xi_n)$ of unit vectors in $L^2(A)$ such that $\lim_n\norm{\pi^{\alpha}_A(\mu)\xi_n - c\xi_n} = 1$
(see the proof of \cite[Th\'eor\`eme 1]{DG}). Using Lemma \ref{positiveL2} and the Cauchy-Schwarz 
inegality, we get
\begin{align*}
\Big|\int_G&\scal{\xi_n,\pi^{\alpha}_A(s)\xi_n}d\mu(s)\Big|\leq \int_G\abs{\scal{\xi_n,\pi^{\alpha}_A(s)\xi_n}}d\mu(s)\\
&\leq \int_G\scal{\abs{\xi_n}, \pi^{\alpha}_A(s)\abs{\xi_n}}^{1/2}\scal{\abs{J\xi_n},\pi^{\alpha}_A(s)\abs{J\xi_n}}^{1/2} d\mu(s)\\
&\leq \Big(\int_G \scal{\abs{\xi_n}, \pi^{\alpha}_A(s)\abs{\xi_n}}d\mu(s)\Big)^{1/2}
\Big(\int_G \scal{\abs{J\xi_n},\pi^{\alpha}_A(s)\abs{J\xi_n}}d\mu(s)\Big)^{1/2}.
\end{align*}
Since
 $$\int_G \scal{\abs{\xi_n}, \pi^{\alpha}_A(s)\abs{\xi_n}}d\mu(s)\leq 1\quad\hbox{and}\quad \int_G \scal{\abs{J\xi_n},\pi^{\alpha}_A(s)\abs{J\xi_n}}d\mu(s)\leq 1$$
 and since $\dst\lim_n\Big|\int_G\scal{\xi_n,\pi^{\alpha}_A(s)\xi_n}d\mu(s)\Big|= 1$, 
 we infer that 
 $$\dst\lim_n \int_G \scal{\abs{\xi_n}, \pi^{\alpha}_A(s)\abs{\xi_n}}d\mu(s) = 1.$$
It follows that there exists a subsequence $(\eta_n)$ of $(\abs{\xi_n})$ such that
$$\lim_n \scal{\eta_n, \pi^{\alpha}_A(s)\eta_n} = 1$$
and therefore
$$\lim_n\norm{\pi^{\alpha}_A(s)\eta_n - \eta_n} = 0$$
for all $s$ in a subset $S$ of $G$ whose complement has $\mu$-measure zero.

Let us denote by $A^c$ the $C^*$-subalgebra of all $x\in A$ such that $s\mapsto \alpha_s(x)$ is norm
continuous, and for $n\in \N$, denote by $\varphi_n$ the state $x\mapsto \scal{\eta_n, x\eta_n}$
defined on $A^c$. Let $\varphi$ be a weak*-limit point of $(\varphi_n)$ in the dual space of $A^c$.
The set $F$ of elements $s\in G$ such that $\varphi\circ\alpha_s = \varphi$ is a closed subgroup
containing $S$. Therefore we have $\mu(G\setminus F) = 0$. It follows that $F=G$ since
$\mu$ is an adapted probability measure. To conclude we use the well-known fact that the existence
of a $G$-invariant state on $A$ is equivalent to the existence of a $G$-invariant state on $A^c)$
(see \cite[Lemme 2.1]{AD82} for instance).
\end{proof}

\begin{rem} As a consequence of the previous theorem, we see that $\iota_G\prec\pi_A^{\alpha}$
  if and only if $\iota_G\prec\pi_A^{\alpha}\otimes\overline{\pi_A^{\alpha}}$. Indeed, whenever this last condition
  holds, the implication (ii) $\Rightarrow$ (i) gives the existence of a state $\varphi$ on $\cB(L^2(A))$ such that $\varphi\circ\Ad \pi_A^{\alpha}(s)$ for every $s\in G$ and therefore the existence of a $G$-invariant state on $A$
  by restriction. Applied to the regular representation $\lambda_G$,  one recovers a result of Fell \cite{Fe1} saying
  that $G$ is amenable whenever $\lambda_G$ weakly contains a finite dimensional representation.
 \end{rem} 

\begin{rem} As said in the introduction, one cannot expect in general that $\pi_B^{\alpha} \prec \pi_A^{\alpha}$
implies that the action is amenable. This is already not true  when $A$ is abelian. In \cite{AD03}
we studied some particular cases where this fact holds in the abelian setting. Let us recall below another important
particular case, due to Connes \cite{Co}, where this fact holds in the noncommutative setting.
Here we take $A = \cB(L^2(B))$ and $G$ is the unitary group $\cU(B)$ of $B$ equipped with the discrete topology. This group  is not
countable, but we observe that the only result used in the sequel, namely theorem \ref{coamen}, does not require any separability assumption. We let $G$ act on 
$(A,B)$ by $\alpha_U = \Ad U$ for $U\in G$. The amenability of the $G$-action on $(A,B)$ is, by definition,
the injectivity of $B$. Recall that an hypertrace for $B$ is a $G$-invariant state on $A$.
\end{rem}

\begin{thm}[\cite{Co}] We keep the above assumptions. The following conditions are equivalent:
\begin{enumerate}
\item[(i)] There exists an hypertrace for $B$;
\item[(ii)] For every finite subsets $\set{U_1,\dots, U_n}$ of  $\cU(B)$ and $\set{c_1,\dots, c_n}$ of $\C$,
we have
$$\abs{\sum_{i=1}^n c_i} \leq \norm{\sum_{i=1}^n c_i U_i \otimes \overline{U}_i}.$$
\item[(iii)] For every finite subset $\set{U_1,\dots, U_n}$ of $\cU(B)$ we have
$$\norm{\sum_{i=1}^n  U_i \otimes \overline{U}_i} = n.$$
\end{enumerate}
Moreover, if $B$ is a factor, these conditions are equivalent to 
\begin{enumerate}
\item[(iv)] $B$ is a finite injective factor.
\end{enumerate}
In particular, when $B$ is a finite factor, we see that $B$ is injective if and only if $\pi_B^{\alpha} \prec \pi_A^{\alpha}$.
\end{thm}

\begin{proof} Assertion (ii) means that $\iota_G \prec \pi_A^{\alpha}$. Therefore (i) $\Rightarrow$ (ii) is a particular case of (i) $\Rightarrow$ (ii) in theorem \ref{coamen}. 

(ii) $\Rightarrow$ (iii) is obvious. Thanks to (iii) $\Rightarrow$ (i) in theorem \ref{coamen}, for every finite
subset $\cF$ of $\cU(B)$, we get a state $\psi_\cF$ on $A$, invariant by $\Ad U$, $U\in \cF$. Taking
a limit point of $(\psi_\cF)$ along the filter of finite subsets of $\cU(B)$, we obtain an hypertrace for $B$.

 (iv) $\Rightarrow$ (i) is obvious. Indeed, if  $\tau$ is the tracial state of $B$ and
 $E$ is a conditional expectation from $\cB(L^2(B))$ onto $B$, then $\tau\circ E$ is an hypertrace.
 
 Let us sketch the proof of (i)  $\Rightarrow$ (iv) whenever $B$ is a factor. Let $\psi$ be an hypertrace
 for $B$. Its restriction to $B$ is a trace, and therefore $B$ is finite. Denote by $\tau$ its trace. For
 $a\in A^+$ and $b\in B$, we set $\psi_a(b) = \psi(ab)$.  Then $\psi_a$ is a positive state on $B$ with
  $\psi_a \leq \norm{a}\tau$. We denote by $E(a)$ the Radon-Nikod\' ym derivative of $\psi_a$ with respect
  to $\tau$. Then it is easy to check that $E$ extends into a conditional expectation from $A$ onto $B$
  ({\it e.g.}, see \cite[Lemma 2.2]{Ha83}).
  
  Assume now that $B$ is a finite factor. If $B$ is injective, we have  $\pi_B^{\alpha} \prec \pi_A^{\alpha}$ by theorem \ref{faibleinc}.
  Conversely, assume that $\pi_B^{\alpha} \prec \pi_A^{\alpha}$. Since
  the tracial state of $B$ is $G$-invariant, we have $\iota_G\leq \pi_B^{\alpha}$ and therefore (ii) holds.
  It follows that $B$ injective by (iv).
  \end{proof}
  
  The following proposition extends the equivalence between (i) and (ii) stated in theorem \ref{coamen}
  
  \begin{prop} Let $\alpha$ be an action of a locally compact group $G$ on $(A,B)$. We assume that $B$ is a finite factor and that $\alpha(G)$ contains the group
  $\set{\Ad U : U\in \cU(B)}$. Then the action is amenable if and only if $\pi_B^{\alpha} \prec \pi_A^{\alpha}$.
  \end{prop}
  
  \begin{proof} Assume that $\pi_B^{\alpha} \prec \pi_A^{\alpha}$. Since $B$ is a finite factor, we have $\iota_G \leq  \pi_B^{\alpha}$
  and therefore $\iota_G \prec \pi_A^{\alpha}$. By theorem \ref{coamen}, there exists a $G$-invariant state  $\psi$ on $A$.
  We have $\psi(ab) = \psi(ba)$ for $a\in A$ and $b\in \cU(B)$. As in the previous theorem, this state gives rise
  to a conditional expectation from $A$ onto $B$, easily seen to be $G$-invariant.
  \end{proof}

\noindent{\bf Acknowledgements:}  We are grateful to Uffe Haagerup for helpful discussions relative to the proofs
of lemma \ref{Kosaki} and proposition \ref{UH}. We would also like to thank Bachir Bekka for his comments about a
preliminary version and Alain Valette for the reference \cite{Ka}.

\bibliographystyle{alpha}

\end{document}